\documentclass{article}
\usepackage{graphicx} 
\usepackage{amsmath}
\usepackage{mathrsfs}
\usepackage{amssymb}
\usepackage{xcolor}
\usepackage{mathtools}
\usepackage{amsthm}

\theoremstyle{plain}
\newtheorem{theorem}{Theorem}[section]
\theoremstyle{plain}
\newtheorem{definition}{Definition}[section]
\theoremstyle{plain}
\newtheorem{proposition}{Proposition}[section]
\theoremstyle{plain}
\newtheorem{lemma}{Lemma}[section]
\theoremstyle{plain}
\newtheorem{remark}{Remark}[section]

\title{Approximation processes by multidimensional Bernstein-type exponential polynomials on the hypercube}
\author{\textit{Laura Angeloni\thanks{Corresponding author}, Danilo Costarelli, Chiara Darielli} \\
\\
Department of Mathematics and Computer Science \\
            University of Perugia\\
       1, Via Vanvitelli, 06123 Perugia, Italy    \\  \\
  {\small {\tt laura.angeloni@unipg.it}} -  {\small {\tt 
 chiara.darielli@studenti.unipg.it}} \\
 {\small {\tt danilo.costarelli@unipg.it}}  
  }
\date{}

\begin{document}

\maketitle

\begin{abstract}
In this paper we introduce a new family of Bernstein-type exponential polynomials on the hypercube $[0,1]^d$ and study their approximation properties. Such operators fix a multidimensional version of the exponential function and its square. In particular, we prove uniform convergence, by means of two different approaches, as well as a quantitative estimate of the order of approximation in terms of the modulus of continuity of the approximated function.

\vskip0.2cm


\vskip0.3cm

\noindent {\bf AMS Subjclass:} 41A25, 41A30    
\end{abstract}

\section{Introduction}
Polynomial approximation is a classical topic of Approximation Theory, arising from the well-known fundamental Weierstrass theorem (\cite{1885}), asserting that any continuous function within a closed and bounded interval can be approximated through polynomials. One of the most famous proofs of such result is given by means of Bernstein operators, a family of polynomials introduced by Bernstein in 1912. 
Such operators, defined as
\begin{align*} \displaystyle %\label{Ber} 
B_nf(x)=B_n(f,x) := \sum_{k=0}^{n} f\left(\frac{k}{n}\right)\binom{n}{k} x^k(1-x)^{n-k}, \hspace{5mm} x \in [0,1],
\end{align*}
for $f\in C([0,1])$, $n \in \mathbb{N}$,
provide a powerful methodology for approximation, and their properties play a crucial role in this context (see, e.g., \cite{1932,Popo,Freud,ConApp,Carothers,Ber}). 

Starting from the above operators, several generalizations of such polynomials in many different directions have been considered in the literature. Among them, we mention the following remarkable examples of linear and nonlinear operators: the sampling-type operators (\cite{bardaro2017generalization,angeloni2019characterization,bardaro2019exponential,angeloni2021approximation,angeloni2021multivariate}), the max-product operators (\cite{bede2016approximation,boccali2024convergence}), the neural network-type operators (\cite{costarelli2022density,kadak2022multivariate}), the Szasz-Mirakjan operators (\cite{acar2017szasz,2acar2017szasz}) and many others.

\vspace{3mm} 
In the direction of the multidimensional generalization, the definition of Bernstein polynomials was extended to functions $f$ defined on the hypercube $Q_d=[0,1]^d$ as follows:
\begin{align} \displaystyle %\label{BnQd}
    \widetilde{B}_nf(\mathbf{x})=\widetilde{B}_n(f, \mathbf{x}):= \sum_{h_1= 0}^n \ldots & \sum_{h_d = 0}^n
     f\left( \frac{h_1}{n}, \dots, \frac{h_d}{n}\right)\binom{n}{h_1} \dots \binom{n}{h_d} \notag \\ &  x_1^{h_1}(1-x_1)^{n-h_1} \dots x_d^{h_d}(1-x_d)^{n-h_d}, \label{multi-bern-classic}
\end{align} 
$\mathbf{x}=(x_1.\ldots,x_d)\in Q_d$. The approximation properties of $(\widetilde{B}_nf)_n$ were initially explored by 
Hildebrandt and Shoenberg (\cite{Hild}), as well as by Butzer (\cite{TwoDim}).

A remarkable development in Approximation Theory is provided by the Korovkin-type theorems, where the concept of polynomial approximation is generalized to sequences of positive linear operators. Indeed, this theory provides more flexible criteria, allowing to consider broader classes of approximable functions. 

Recently, the study of approximation results by means of exponential-type operators has been deeply developed by several authors (see, e.g., \cite{AKR, acar2017szasz, DAG, GuTa, BerType, AIR, G, Acu, Angeloni}). In particular, in \cite{Morigi} Morigi and Neamtu introduced and studied the following exponential generalization of the Bernstein polynomials: 
\begin{align*} \displaystyle %\label{Gn} 
    \mathscr{G}_nf(x) = \mathscr{G}_n(f,x):= \sum_{k=0}^n f\left( \frac{k}{n} \right) e^{-\mu k / n}e^{\mu x} p_{n,k}(a_n(x)), \ \ n\in \mathbb{N},
\end{align*} 
for $f\in C([0,1])$, where
\begin{align*} %\label{an} 
\displaystyle
    a_n(x)=\frac{e^{\mu x/n}-1}{e^{\mu/n}-1},
\end{align*}
$\mu >0$, and $p_{n,k}(x)=\binom{n}{k} x^k (1-x)^{n-k}$. Besides the uniform convergence, they proved shape preserving properties and a Voronovskaja formula. 

Following the above mentioned path for the classical Bernstein polynomials, we introduce in this paper a multidimensional version of the operators $G_n$ for functions defined on the hypercube $Q_d$ and study their approximation properties. In particular, for such operators we prove the uniform convergence in case of continuous functions by means of the concept of Korovkin subset. Moreover, we also furnish a proof of the same result by means of a constructive approach, since this opens the way to the achievement of a quantitative estimate employing the modulus of continuity of the approximated function. As a direct consequence it is possible to obtain the corresponding qualitative order of approximation in case of functions belonging to suitable Lipschitz classes.

\section{Preliminaries}
Let us provide some notations, definitions and properties that will be used throughout the paper.

Given a metric space $(X,d)$, we denote by $F(X)$ the linear space of all the real-valued functions defined on $X$ and by $C(X)$ the subspace of $F(X)$ of all the continuous functions on $X$.

\begin{definition} 
A linear subspace $E$ of $F(X)$ is said to be a lattice subspace if
\begin{align*} \displaystyle 
|f|\in E \hspace{4mm} for \hspace{2mm} every \hspace{2mm} f \in E. 
\end{align*} 
\end{definition}

\begin{definition}
    A Banach lattice $E$ is a vector space endowed with a norm $\Vert \cdot \Vert$ and an order relation $\leq$ on $E$ such that
    \begin{itemize}
        \item[(i)] $\left(E,\Vert \cdot \Vert\right)$ is a Banach space;
        \item[(ii)] $\left(E,\leq \right)$ is a vector lattice;
        \item[(iii)] if $f,g \in E$ and $\vert f \vert \leq \vert g \vert$, then $\Vert f \Vert \leq \Vert g \Vert$. \end{itemize} 
\end{definition}

\begin{definition}
    Given a metric space $(Y,d')$, a linear operator $T: E \rightarrow F(Y)$ is said to be positive if \begin{align*} \displaystyle
    T(f) \ge 0 \hspace{4mm} for \hspace{2mm} every \hspace{2mm} f \in E, \hspace{2mm} f \ge 0.
\end{align*}
\end{definition}

We recall the celebrated Korovkin theorem, that provides a very useful and simple criterion in order to prove that a given sequence $(T_n)_{n\ge1}$ of positive linear operators on $C([0,1])$ is an approximation process, i.e., $T_n(f) \rightarrow f$ uniformly on $[0,1]$, as $n\rightarrow +\infty$, for every $f \in C([0,1])$.

\begin{theorem}[Korovkin, 1953] \label{K1} Let $(T_n)_{n\ge1}$ be a sequence of positive linear operators from $C([0,1])$ into $C([0,1])$ such that, for every $g \in \{e_0, e_1, e_2\}$, with $e_i(t)=t^i$, $i=0,1,2,$
\begin{align*} \displaystyle
    \lim_{n \rightarrow +\infty} T_n(g) = g \hspace{5mm} uniformly \hspace{3mm} on \hspace{3mm} [0,1].
\end{align*} Then, for every $f \in C([0,1])$, 
\begin{align*} \displaystyle
    \lim_{n \rightarrow +\infty} T_n(f) = f \hspace{5mm} uniformly \hspace{3mm} on \hspace{3mm} [0,1].
\end{align*} 
\end{theorem}

In \cite{Bohman}, H. Bohman showed a result analogous to Theorem \ref{K1} by considering sequences of positive linear operators on $C([0,1])$ of the form 
\begin{align*} \displaystyle
    T(f,x) = \sum_{i \in I} f(a_i)\varphi_i(x), \hspace{4mm} 0 \leq x \leq 1,
\end{align*} where $(a_i)_{i \in I}$ is a finite family in $[0,1]$, $\varphi_i \in C([0,1])$ and $I$ is a set of indices.
For this reason, Theorem \ref{K1} is often called the Bohman-Korovkin theorem.

The above theorem can be generalized introducing the notion of Korovkin subset (see, e.g., \cite{AltoCam,Korovkin}).

\begin{definition}
    Let $E$ and $F$ be Banach lattices and consider a positive linear operator $T: E \rightarrow F.$ A subset $M$ of $E$ is said to be a {\rm Korovkin subset} of $E$ for $T$ if for every sequence $(T_n)_{n\ge 1}$ of positive linear operators from $E$ into $F$ satisfying 
    \begin{itemize}
      \item[(i)] $\sup\limits_{n \ge 1} \Vert T_n \Vert < +\infty$
    \end{itemize} and
    \begin{itemize}
         \item[(ii)] $\lim\limits_{n\rightarrow+\infty} T_n(g)=T(g)$ for every $g \in M$,
    \end{itemize} it turns out that 
    \begin{align*}
        \lim_{n \rightarrow +\infty} T_n(f)=T(f) \hspace{3mm} for \hspace{2mm} every \hspace{2mm} f\in E.
    \end{align*}
\end{definition}

By means of such definition, as a consequence of Theorem \ref{K1}, we can assert that $e_i(t)=t^i$, $i=0,1,2$ form a Korovkin subset of $C([0,1])$. A useful criterion to find Korovkin subsets is furnished by the following proposition.

\begin{proposition} \label{Prop64}
    Let $M$ be a subset of $C(X)$ and let $\mathcal{L}(M)$ be the linear subspace of $C(X)$ generated by $M$. Assume that for every $x,y \in X$, $x \neq y$, there exists $h \in \mathcal{L}(M)$, $h \ge 0$, such that $h(x)=0$ and $h(y)>0$. Then $M$ is a Korovkin subset of $C(X)$.
\end{proposition}

\begin{remark} \label{remark} 
\textup{Proposition \ref{Prop64} also holds for any compact metric space $X\subseteq \mathbb{R}^d$.}
\end{remark}

Applying the previous proposition it is easy to find other Korovkin subsets of $C([0,1])$, besides the power functions $\{1,x,x^2\}$, such as, for example, $\{1,\exp({\lambda_1 x}),\exp({\lambda_2 x})\}$, $x \in [0,1]$, where $0<\lambda_1<\lambda_2$.

In case of functions of several variables, it is possible to prove that 
$$
\left\{e_0,pr_1,\dots,pr_d,\sum_{i=1}^d pr_i^2 \right\},
$$
where $e_0({\bf x}):=1$ and $pr_i({\bf x}):=x_i$, for every ${\bf x}:=(x_1,\ldots,x_d)\in X$, $X$ compact subset of $\mathbb{R}^d$, $i=1,\ldots,d$, is a Korovkin subset of $C([0,1]^d)$ (see \cite{Korovkin}).

For details and insights on Korovkin-type approximation theory, see, e.g., \cite{Bohman,1953,AltoCam,Korovkin} and the references therein.

\section{Bernstein-type Exponential Polynomials}

In \cite{Morigi} the authors introduced a generalization of the Bernstein operators associated to an exponential function. For a fixed real parameter $\mu > 0$, by $\exp_\mu$ we denote $\exp_\mu(x)=e^{\mu x}$, $x \in \mathbb{R}$.
\vspace{2mm}\\ 
The above mentioned generalization is defined for $f\in C([0,1])$ and for $n \in \mathbb{N}$ by
\begin{align} \label{Gn} \displaystyle
    \mathscr{G}_nf(x) = \mathscr{G}_n(f,x):= \sum_{k=0}^n f\left( \frac{k}{n} \right) e^{-\mu k / n}e^{\mu x} p_{n,k}(a_n(x)),
\end{align} where
\begin{align*} %\label{an} 
\displaystyle
    a_n(x)=\frac{e^{\mu x/n}-1}{e^{\mu/n}-1}.
\end{align*}
The connection between the operators (\ref{Gn}), called Bernstein-type exponential operators, and the classical Bernstein operators is given by
\begin{align} \label{Rel} \displaystyle
    \mathscr{G}_n (f,x) = \exp_\mu (x) B_n \left( \frac{f}{\exp_\mu}, a_n(x)\right).
\end{align}
We note that, for each $n \in \mathbb{N}$, $a_n(x)$ is an increasing and convex continuous function that satisfies $a_n(0)=0, \hspace{1mm} a_n(1)=1$ and $a_n(x) >0$ for every $x \in ]0,1]$. Consequently, $\mathscr{G}_n$ is a positive operator that interpolates continuous functions at the endpoints of the interval $[0,1]$. Moreover we note that the operators $\mathscr{G}_n$ reproduce $\exp_\mu$ and $\exp^2_\mu$, i.e.,
\begin{align} \label{e12} \displaystyle
    \mathscr{G}_n(\exp_\mu,x)=e^{\mu x}, \hspace{10mm} \mathscr{G}_n(\exp^2_\mu,x)=e^{2\mu x}.
\end{align}
We point out that the second relation is a consequence of the equality 
$$
\sum_{k=0}^n e^{\mu k\over n} p_{n,k}(a_{n}(x))=e^{\mu x }
$$ (see \cite{AOR}).

In addition, we recall the following lemma showing some basic identities of the operators $\mathscr{G}_n$ that will be used later. %They can be derived by direct calculations, or with the use of some mathematical software. 
\begin{lemma}[\cite{BerType}]\label{Lemma5} 
For each $n \in \mathbb{N}$ and $x \in [0,1]$, the following identities hold:
\begin{itemize}
    \item[(i)] $\mathscr{G}_n(e_0,x) = e^{\mu (x-1)} \left(e^{\mu/n}+1-e^{\mu x/n}\right)^n$,
    \item[(ii)] $\mathscr{G}_n(\exp^3_\mu,x) = e^{\mu x} \left(e^{\mu(x+1)/n}+e^{\mu x/n}-e^{\mu/n}\right)^n$,
    \item[(iii)] $\mathscr{G}_n(\exp^4_\mu,x) = e^{\mu x} \left(e^{\mu(x+2)/n}+ e^{\mu(x+1)/n}+e^{\mu x/n}-e^{\mu/n}-e^{2\mu/n}\right)^n.$
\end{itemize} \end{lemma} 
For each $x \in (0,1)$, we consider the function $\exp_{\mu,x}$ defined for every $t \in [0,1]$ by 
\begin{align*} \displaystyle \displaystyle
 \exp_{\mu,x}(t):=e^{\mu t}-e^{\mu x}.
\end{align*}
%By elementary calculation, we find that, for each $x \in (0,1)$ and $t \in [0,1]$, whenever $\mu \ge 1$, % ogni volta che 
%\begin{align*} \displaystyle
%    e_x^2(t) \leq \exp_{\mu,x}^2(t),
%\end{align*} where $e_x^2(t) = (t-x)^2$ and $\exp_{\mu,x}^2(t)=(e^{\mu t}-e^{\mu x})^2$. 
Using the expression in (\ref{Rel}) and Lemma \ref{Lemma5}, we can explicitly compute the operators $\mathscr{G}_n$ in $\exp_{\mu,x}^2$, for every $x \in [0,1]$:
\begin{align} \label{exp2} \displaystyle
    \mathscr{G}_n(\exp^2_{\mu,x},x) & %= \mathscr{G}_n(\exp^2_\mu,x)-2e^{\mu x}\mathscr{G}_n(\exp_\mu,x)+e^{2\mu x}\mathscr{G}_n (e_o,x)\notag \\ & =e^{2\mu x}-2e^{2\mu x}+e^{2\mu x}\mathscr{G}_n (e_o,x)\notag \\ & 
    =e^{2\mu x}\mathscr{G}_n (e_o,x)-e^{2\mu x} = e^{2\mu x}(\mathscr{G}_n (e_o,x)-1) 
    \notag \\ & = e^{2\mu x}\left( e^{\mu(x-1)}\left( e^{\mu/n} +1- e^{\mu x/n} \right)^n -1 \right).
\end{align}

In \cite{BerType} the authors prove the following theorem of uniform convergence for the operators $(\mathscr{G}_n f)_n$. 
\begin{theorem} \label{Gnf-f}
    If $f \in C([0,1])$, then $\mathscr{G}_n f$ converges to $f$ uniformly on $[0,1]$ as $n\rightarrow +\infty$.
\end{theorem}
Furthermore, also a Voronovskaja formula has been proved. In particular, in the case of $e_0$, one has:
\begin{align*} %\label{x(1-x)} 
\displaystyle
    \lim_{n \rightarrow \infty} n(\mathscr{G}_n (e_0,x)-1) & = \lim_{n \rightarrow \infty} n \left( e^{\mu(x-1)} \left(e^{\mu/n}+1-e^{\mu x/n}  \right)^n-1\right)
    \notag \\ & = \mu^2x(1-x), \ x\in [0,1]. 
\end{align*}
In particular, from the above limit it is easy to deduce the following useful (and uniform) inequality:
\begin{equation} \label{voro-e0}
    |\mathscr{G}_n (e_0,x)-1|\ \leq {\mu^2 \over n}, \quad x \in [0,1], \quad n \in \mathbb{N}.
\end{equation}
For additional information regarding shape preserving characteristics, an analysis of Bernstein polynomials applied to convex functions and a study of the problem of simultaneous approximation see \cite{Arama,Mond,schoe,Convexity,BerType,Acu}.

\section{Multidimensional Bernstein-type Exponential Polynomials}

We now extend the definition of the operators $\mathscr{G}_n$ to obtain a new tool to reconstruct multivariate functions. The idea is to define the operators $\mathscr{G}_n$ in the multidimensional case, with the purpose of studying some approximation properties in this setting. 

We consider the hypercube $Q_d=[0,1]^d$ and the vector $\mathbf{x}=(x_1,\dots,x_d) \in \mathbb{R}^d$, where $d \ge 1$. Additionally, we let 
\begin{align} \label{expmud}
    \exp_\mu(\mathbf{x}):= e^{\mu x_1} \dots \hspace{1mm} e^{\mu x_d}= e^{\mu \sum_{i=1}^d x_i},
\end{align} 
that is, the multidimensional exponential function defined on the space $\mathbb{R}^d$.

For every $n \ge 1$, $f\in C(Q_d)$ and $\mathbf{x} \in Q_d$, we can define the multidimensional Bernstein-type exponential operators as follows:
\begin{align} \label{GnD} \displaystyle
\widetilde{\mathscr{G}}_nf(\mathbf{x})=\widetilde{\mathscr{G}}_n(f, \mathbf{x}):= {\sum_{k_1=0}^n} & \dots {\sum_{k_d=0}^n}
\hspace{1mm} f\left( \frac{k_1}{n}, \dots, \frac{k_d}{n}\right)
e^{-\mu k_1/n} \dots \hspace{1mm} e^{-\mu k_d/n} \notag \\ & e^{\mu x_1} \dots \hspace{1mm} e^{\mu x_d} \hspace{1mm} p_{n,k_1}(a_{n}(x_1)) \dots \hspace{1mm}p_{n,k_d}(a_{n}(x_d)).
\end{align}
The definition of the multidimensional Bernstein polynomial in the variables $x_1,\dots,x_d$, corresponding to the function $f\in C(Q_d)$, is provided in (\ref{multi-bern-classic}). 

Even in the multidimensional case, there exists a connection between the operators defined in (\ref{GnD}) and those ones defined in (\ref{multi-bern-classic}), expressed as
\begin{align} \label{Conn} \displaystyle
 &   \widetilde{\mathscr{G}}_n(f, \mathbf{x})  = {\sum_{k_1=0}^n} \dots {\sum_{k_d=0}^n} \hspace{1mm}
    f\left( \frac{k_1}{n}, \dots, \frac{k_d}{n}\right)
    e^{-\mu k_1/n} \dots \hspace{1mm} e^{-\mu k_d/n} \notag \\
    & \hspace{14mm} e^{\mu x_1} \dots \hspace{1mm} e^{\mu x_d} \hspace{1mm} p_{n,k_1}(a_{n}(x_1)) \dots \hspace{1mm} p_{n,k_d}(a_{n}(x_d)) \notag \\
    & =: \hspace{1mm} e^{\mu x_1} \dots \hspace{1mm} e^{\mu x_d} {\sum_{k_1=0}^n} \dots {\sum_{k_d=0}^n} \hspace{1mm} f_\mu \left( \frac{k_1}{n}, \dots, \frac{k_d}{n}\right) p_{n,k_1}(a_{n}(x_1)) \dots \hspace{1mm} p_{n,k_d}(a_{n}(x_d)) \notag \\ & = \hspace{1mm} e^{\mu x_1} \dots \hspace{1mm} e^{\mu x_d} \hspace{1mm} \widetilde {B}_n f_\mu(a_n(x_1),\dots,a_n(x_d)).
\end{align} 
where
\begin{align*} \displaystyle
    f_\mu\left( \frac{k_1}{n}, \dots, \frac{k_d}{n}\right) := f\left( \frac{k_1}{n}, \dots, \frac{k_d}{n}\right) e^{-\mu k_1/n} \dots \hspace{1mm} e^{-\mu k_d/n}. 
    %= e^{-\mu \sum_{i=1}^d k_i/n} = \exp_\mu ^{-1}\left(\frac{k_1}{n},\dots ,\frac{k_d}{n}\right).
\end{align*} 
While $B_n$ fixes the functions $e_0$ and $e_1$, the operator $\mathscr{G}_n$ reproduces $\exp_\mu$ and $\exp^2_\mu$ and this property holds true in the multidimensional case as well. Indeed, the operator $\widetilde{\mathscr{G}}_n$, defined in (\ref{GnD}), fixes the multidimensional exponential functions $\exp_\mu$ and $\exp_\mu^2$. In particular, as in (\ref{Conn}), we have that
\begin{align*} %\label{rel1D} 
\displaystyle
    \widetilde{\mathscr{G}}_n(\exp_\mu, \mathbf{x}) & = \exp_\mu(\mathbf{x}) \sum_{k_1=0}^n {\dots} \sum_{k_d=0}^n \hspace{1mm} p_{n,k_1}(a_n(x_1)) \hspace{1mm} \dots \hspace{1mm} p_{n,k_d}(a_n(x_d)) \notag \\ & = \exp_\mu(\mathbf{x}) \left( \hspace{1mm}\sum_{k_1=0}^n \hspace{1mm} p_{n,k_1}(a_n(x_1)) \hspace{1mm} {\dots} \sum_{k_d=0}^n p_{n,k_d}(a_n(x_d))\right) \notag \\
    & = \exp_\mu(\mathbf{x}),
    \end{align*} 
    taking into account that $\displaystyle\sum_{k=0}^n p_{n,k}(a_n(x))=1$, for every $x \in [0,1]$.  For the other multidimensional exponential function, there holds 
\begin{align} \label{rel2D} \displaystyle
    \widetilde{\mathscr{G}}_n(\exp_\mu^2, \mathbf{x}) & = {\sum_{k_1=0}^n} \dots {\sum_{k_d=0}^n} \exp_\mu^2 \left( \frac{k_1}{n},{\dots},\frac{k_d}{n}\right) e^{-\mu k_1/n} \hspace{1mm} \dots \hspace{1mm} e^{-\mu k_d /n}\notag \\ 
    & \hspace{14mm} e^{\mu x_1} \hspace{1mm} \dots \hspace{1mm} e^{\mu x_d} \hspace{1mm} p_{n,k_1}(a_{n}(x_1)) \hspace{1mm} \dots\hspace{1mm} p_{n,k_d}(a_{n}(x_d)) \notag \\
    & = {\sum_{k_1=0}^n} \dots {\sum_{k_d=0}^n} \hspace{1mm} e^{2\mu k_1 /n} \hspace{1mm} \dots \hspace{1mm} e^{2\mu k_d/n} \hspace{1mm} e^{-\mu k_1/n} \hspace{1mm} \dots \hspace{1mm} e^{-\mu k_d /n} \notag \\ 
    & \hspace{14mm} e^{\mu x_1} \hspace{1mm} \dots \hspace{1mm} e^{\mu x_d} \hspace{1mm} p_{n,k_1}(a_{n}(x_1)) \hspace{1mm} \dots \hspace{1mm} p_{n,k_d}(a_{n}(x_d)) \notag \\
    & = \sum_{k_1=0}^n e^{\mu k_1 /n}\hspace{1mm}e^{\mu x_1}\hspace{1mm} p_{n,k_1}(a_{n}(x_1)) \hspace{1mm} \dots \hspace{1mm} \sum_{k_d=0}^n e^{\mu k_d /n}\hspace{1mm}e^{\mu x_d}\hspace{1mm} p_{n,k_d}(a_{n}(x_d)) \notag \\[0.5em]
    & = {\mathscr{G}}_n(\exp_\mu^2, x_1) \hspace{1mm} \dots \hspace{1mm}{\mathscr{G}}_n(\exp_\mu^2, x_d) = \exp_{\mu}^2(x_1) \hspace{1mm} \dots \hspace{1mm} \exp_{\mu}^2(x_d) \notag \\[0.5em] & =\exp_\mu^2(\mathbf{x}),
\end{align} by (\ref{e12}) and (\ref{expmud}). 
In conclusion, these operators fix the multidimensional exponential functions, i.e., $\widetilde{\mathscr{G}}_n(\exp_\mu, \mathbf{x})=\exp_\mu(\mathbf{x})$ and $\widetilde{\mathscr{G}}_n(\exp_\mu^2, \mathbf{x})=\exp_\mu^2(\mathbf{x})$. 
\vspace{2mm} \\ 
It is clear that, similarly to (\ref{rel2D}), the operators $\widetilde{\mathscr{G}}_n$ applied to $\exp^m_\mu$ reduce to 
\begin{align} \label{produt}
    \widetilde{\mathscr{G}}_n(\exp_\mu^m,\mathbf{x})= \prod_{i=1}^d \mathscr{G}_n(\exp_\mu^m,x_i),
\end{align} for every positive integer $m$ and ${\mathbf x}\in Q_d$. Then the following result can be easily derived by Lemma \ref{Lemma5}.
\begin{lemma} \label{LemD}
    For $\mathbf{x} \in Q_d$, $n \in \mathbb{N}$, there holds 
    \begin{itemize}
    \item[(a)] $\widetilde{\mathscr{G}}_n(e_0, \mathbf{x})=e^{\mu(x_1-1)}(e^{\mu/n}+1-e^{\mu x_1/n})^n \hspace{1mm} \dots \hspace{1mm} e^{\mu (x_d-1)}(e^{\mu/n}+1-e^{\mu x_d/n})^n,$
    \item[(b)] $\widetilde{\mathscr{G}}_n(\exp_\mu^3, \mathbf{x})=e^{\mu x_1} \left(e^{\mu(x_1+1)/n}+e^{\mu x_1/n}-e^{\mu/n}\right)^n \hspace{1mm} \dots $
    
    $\hspace{47mm} e^{\mu x_d} \left(e^{\mu(x_d+1)/n}+e^{\mu x_d/n}-e^{\mu/n}\right)^n,$
    \item[(c)] $\widetilde{\mathscr{G}}_n(\exp_\mu^4, \mathbf{x})= e^{\mu x_1}\left(e^{\mu (x_1+2)/n}+e^{\mu(x_1+1)/n}+e^{\mu x_1/n}-e^{\mu/n}-e^{2\mu/n}\right)^n \vspace{2mm} \\[0.5em] \dots \hspace{1mm} e^{\mu x_d}\left(e^{\mu (x_d+2)/n}+e^{\mu(x_d+1)/n}+e^{\mu x_d/n}-e^{\mu/n}-e^{2\mu/n}\right)^n.$
    \end{itemize}
\end{lemma}

\section{Approximation Properties}

We now discuss and prove convergence results for the operators $\widetilde{\mathscr{G}}_n$. We first study the uniform convergence of our operators. In order to do it we preliminarily prove the uniform convergence for the function $e_0$. 
\begin{theorem} \label{teoe0}
    $\widetilde{\mathscr{G}}_n e_0$ converges to $e_0$ uniformly on $Q_d$, as $n \rightarrow +\infty$.
\end{theorem}
\begin{proof}[Proof.]
From the equation (\ref{produt}), we obtain
\begin{align*} \displaystyle
    \big|\widetilde{\mathscr{G}}_n(e_0,\mathbf{x}) - e_0(\mathbf{x})\big| & \leq 
     \prod_{i=2}^d \hspace{1mm} \mathscr{G}_n(e_0,x_i) \hspace{1mm} \big| \mathscr{G}_n(e_0,x_{1}) - 1 \big| + \prod_{i=3}^d \hspace{1mm} \mathscr{G}_n(e_0,x_i) \hspace{1mm} \big| \mathscr{G}_n(e_0,x_{2}) - 1 \big| \\[0.5em]
     & \hspace{5mm} + \hspace{1mm} {\dots} \hspace{1mm} + \hspace{1mm} \mathscr{G}_n(e_0,x_d) \hspace{1mm} \big| \mathscr{G}_n(e_0,x_{d-1}) - 1 \big| + \big| \mathscr{G}_n(e_0,x_{d}) - 1 \big| \\[0.5em]
     & = \sum_{j=2}^d \left[ \hspace{1mm} \prod_{i=j}^d \hspace{1mm}\mathscr{G}_n(e_0,x_i) \hspace{1mm} \big| \mathscr{G}_n(e_0,x_{j-1}) - 1 \big| \right] + \big|\mathscr{G}_n(e_0,x_{d}) - 1 \big|.
\end{align*} 
By $(i)$ of Lemma \ref{Lemma5}, we have that 
\begin{align*} %\label{eallamu} 
\displaystyle
    \big| \mathscr{G}_n(e_0,x_i) \big| \leq e^\mu,
\end{align*}
for every $x_i \in [0,1]$, with $i=1,\dots,d$, and for every $n \in \mathbb{N}$. 
In addition, we note that, recalling the consequence of the Voronovskaja formula stated in (\ref{voro-e0}), 
we finally obtain that
\begin{align*} \displaystyle
    \big|\widetilde{\mathscr{G}}_n(e_0,\mathbf{x}) - e_0(\mathbf{x})\big| & \leq \hspace{1mm} e^{\mu(d-1)} \hspace{1mm}\big| \mathscr{G}_n(e_0,x_{1}) - 1 \big| + e^{\mu(d-2)} \hspace{1mm}\big| \mathscr{G}_n(e_0,x_{2}) - 1 \big| \\[0.5em]
    & \hspace{5mm} + \hspace{1mm} \dots \hspace{1mm} + e^{\mu} \hspace{1mm}\big| \mathscr{G}_n(e_0,x_{d-1}) - 1 \big| + \big| \mathscr{G}_n(e_0,x_{d}) - 1 \big| \\[0.3em]
    & \leq \hspace{1mm}\frac{\mu^2}{n} \hspace{1mm} \sum_{j=0}^{d-1} e^{\mu j},
\end{align*} 
for $n$ large enough and this completes the proof. \\ 
\end{proof}
Now, we can establish the following result thanks to the use of a Korovkin theorem.
\begin{theorem} 
    If $f \in C(Q_d)$, then $\widetilde{\mathscr{G}}_n f$ converges to $f$ uniformly on $Q_d$, as $n \rightarrow +\infty$.
\end{theorem}
\begin{proof}
    We notice that the set $M:=\{e_0,pr_{\mu,1},\dots,pr_{\mu,d},pr^2_{\mu,1},\dots,pr^2_{\mu,d}\}$, where 
\begin{align*} %\label{pri}
    pr_{\mu,i}(\mathbf{x}):=\exp_\mu(x_i), \hspace{8mm} i=1,\dots,d, 
\end{align*} 
    with $\mathbf{x} \in Q_d$, is a Korovkin subset of $C(Q_d)$, by Proposition \ref{Prop64}. Indeed, if
    we fix $\widetilde{\bf{x}} \in Q_d$, 
    %then with similar reasonings to Proposition \ref{PropNuova} (one-dimensional case) 
    it is easy to see that the function
    \begin{align*}
        h(\mathbf{x}):= d-2{\sum_{i=1}^d e^{\mu(x_i-\widetilde{x}_i)}}
    +\sum_{i=1}^d e^{2\mu(x_i-\widetilde{x}_i)}    
    \end{align*}
    is such that $h(\mathbf{x})\ge 0$, for every $\mathbf{x} \in Q_d$, and $h(\mathbf{\widetilde{x}})=0,$ while $h(\mathbf{x})>0,$ for every $\mathbf{x}\neq \mathbf{\widetilde{x}}$. Therefore, the thesis is fulfilled if we check the uniform convergence for the functions of $M$. 
    \vspace{2mm} \\ 
    The uniform convergence on $e_0$ is given by  Theorem \ref{teoe0}. Instead, for the functions $pr_{\mu,j}$, with $j=1,\dots,d$, we note that
\begin{align*}
  &  \widetilde{\mathscr{G}}_n (pr_{\mu,j},\mathbf{x}) = \sum_{k_1=0}^n \hspace{1mm}\dots \hspace{1mm} \sum_{k_d=0}^n \hspace{1mm} e^{\mu{k_j}/{n}}\hspace{1mm}e^{-\mu{k_1}/{n}}\dots e^{-\mu k_d/n} \\[0.5em] & \hspace{14mm}
        e^{\mu x_1}\dots e^{\mu x_d}\hspace{1mm}p_{n,k_1}(a_n(x_1))\dots p_{n,k_d}(a_n(x_d)) \\[0.5em] 
    & = \sum_{k_j=0}^n \hspace{1mm} e^{\mu k_j/n}e^{-\mu k_j/n}\hspace{1mm} e^{\mu x_j} p_{n,k_j}(a_n(x_j)) \hspace{1mm} \left[ \prod_{\substack{i=1 \\ i\neq j}}^d e^{\mu x_i} \sum_{k_i=0}^n e^{-\mu k_i/n}p_{n,k_i}(a_n(x_i))\right] \\[0.5em]
    &  = \mathscr{G}_n(\exp_\mu,x_j)\hspace{1mm}\prod_{\substack{i=1 \\ i\neq j}}^d \mathscr{G}_n(e_0,x_i) = \exp_\mu(x_j)\hspace{1mm}\prod_{\substack{i=1 \\ 
    i\neq j}}^d \mathscr{G}_n(e_0,x_i),
\end{align*} 
taking into account that $\mathscr{G}_n(\exp_\mu, x_j)= \exp_\mu(x_j)$ from the first equation of (\ref{e12}). Hence, we obtain 
\begin{align*}
    \widetilde{\mathscr{G}}_n(pr_{\mu,j},\mathbf{x})=\exp_{\mu}(x_j)\prod_{\substack{i=1 \\ i\neq j}}^d \mathscr{G}_n(e_0,x_i), 
\end{align*} that tends to $\exp_{\mu}(x_j)=pr_{\mu,j}(\mathbf{x})$, as $n \rightarrow+\infty$, since Theorem \ref{Gnf-f} holds. \vspace{2mm} \\ Now, for the remaining functions $pr^2_{\mu,j}$, with $j=1,\dots,d$, we do similar reasonings, i.e.,
\begin{align*} 
\widetilde{\mathscr{G}}_n (pr^2_{\mu,j},\mathbf{x}) &= \sum_{k_1=0}^n \hspace{1mm}\dots \hspace{1mm} \sum_{k_d=0}^n \hspace{1mm} e^{2\mu{k_j}/{n}}\hspace{1mm}e^{-\mu{k_1}/{n}}\dots e^{-\mu k_d/n} \\[0.5em] & \hspace{14mm}
        e^{\mu x_1}\dots e^{\mu x_d}\hspace{1mm}p_{n,k_1}(a_n(x_1))\dots p_{n,k_d}(a_n(x_d)) \\[0.5em] 
    & = \sum_{k_j=0}^n \hspace{1mm} e^{2\mu k_j/n}e^{-\mu k_j/n} e^{\mu x_j}\hspace{1mm} p_{n,k_j}(a_n(x_j)) \hspace{1mm} \prod_{\substack{i=1 \\ i\neq j}}^d e^{\mu x_i} \sum_{k_i=0}^n e^{-\mu k_i/n}p_{n,k_i}(a_n(x_i)),
\end{align*} where the first sum is equal to $\mathscr{G}_n(\exp^2_\mu,x_j)= \exp_\mu^2(x_j)$ from the second equation of (\ref{e12}). So we obtain
\begin{align*}
\widetilde{\mathscr{G}}_n(pr^2_{\mu,j},\mathbf{x})=\exp_{\mu}^2(x_j)\prod_{\substack{i=1 \\ i\neq j}}^d \mathscr{G}_n(e_0,x_i), 
\end{align*}
that tends to $\exp_\mu^2(x_j)=pr^2_{\mu,j}(\mathbf{x})$, as $n \rightarrow +\infty$, from Theorem \ref{Gnf-f}. This completes the proof. \\
\end{proof}
Now, we will provide an alternative proof of the previous result using a different technique, that is, my means of a constructive approach. This will be the starting point to achieve quantitative estimates for the order of approximation. % cioè un approccio costruttivo
\begin{proof}
Without any loss of generality, we will prove the result in the bi-dimensional case, that is, $d=2$, for the sake of clearness. % per motivi di chiarezza
Of course, the general case is completely analogous. 
\vspace{2mm} \\
Taking $\mathbf{x} = (x_1,x_2)$, we get
\begin{align*}
    \big|\widetilde{\mathscr{G}}_n(f,\mathbf{x}) - f(\mathbf{x})\big| &\leq
    \big|\widetilde{\mathscr{G}}_n(f,\mathbf{x}) - f(\mathbf{x})\hspace{1mm}\widetilde{\mathscr{G}}_n(e_0,\mathbf{x}) \big| + \big|f(\mathbf{x})\hspace{1mm}\widetilde{\mathscr{G}}_n(e_0,\mathbf{x}) - f(\mathbf{x})\big| \\[0.5em] & =: I_1+I_2.
\end{align*} Now let us fix $\varepsilon > 0$ and we denote by $\gamma >0$ the parameter of the uniform convergence of $f$, i.e.,
\begin{align*}
    \left|f(x_1,x_2)-f(y_1,y_2) \right|\leq \varepsilon, \hspace{5mm} \forall \hspace{1mm} (x_1,x_2),(y_1,y_2) \in Q_d,
\end{align*} 
with $\Vert (x_1,x_2)-(y_1,y_2)\Vert_2 \leq \gamma$, where $\Vert \cdot \Vert_2$ denotes the usual Euclidean norm. Then we have 
\begin{align*}
    I_1 &\leq \sum_{k_1=0}^n \sum_{k_2=0}^n \left| f\left( \frac{k_1}{n},\frac{k_2}{n}\right)-f(x_1,x_2)\right| e^{-\mu k_1/n}e^{-\mu k_2/n} \\[0.5em] & \hspace{24mm} e^{\mu x_1} e^{\mu x_2} \hspace{1mm} p_{n,k_1}(a_n(x_1))p_{n,k_2}(a_n(x_2)) \\[0.5em]
    & = \Bigg( \sum_{\left|\frac{k_1}{n}-x_1 \right|\leq \frac{\gamma}{\sqrt{2}}} \hspace{1mm}\sum_{\left|\frac{k_2}{n}-x_2\right|\leq \frac{\gamma}{\sqrt{2}}} + \sum_{\left|\frac{k_1}{n}-x_1 \right|\leq \frac{\gamma}{\sqrt{2}}} \hspace{1mm}\sum_{\left|\frac{k_2}{n}-x_2\right|> \frac{\gamma}{\sqrt{2}}} +
     \\[0.5em]
     & \hspace{10mm} + \sum_{\left|\frac{k_1}{n}-x_1 \right|> \frac{\gamma}{\sqrt{2}}} \hspace{1mm}\sum_{\left|\frac{k_2}{n}-x_2\right|\leq \frac{\gamma}{\sqrt{2}}} + \sum_{\left|\frac{k_1}{n}-x_1 \right|> \frac{\gamma}{\sqrt{2}}} \hspace{1mm}\sum_{\left|\frac{k_2}{n}-x_2\right|> \frac{\gamma}{\sqrt{2}}} \Bigg) \\[1em]
    & \hspace{8mm} \left| f\left( \frac{k_1}{n},\frac{k_2}{n}\right)-f(x_1,x_2)\right| e^{-\mu k_1/n}e^{-\mu k_2/n}e^{\mu x_1} e^{\mu x_2}  \\[0.5em] & \hspace{9mm} p_{n,k_1}(a_n(x_1))\hspace{1mm}p_{n,k_2}(a_n(x_2)) =: I_{1,1}+I_{1,2}+I_{1,3}+I_{1,4}.
\end{align*} 
We note that
\begin{align*}
    I_{1,1} \leq \varepsilon \hspace{1mm} \left(e^{\mu} \right)^2 \hspace{1mm} \sum_{k_1=0}^n p_{n,k_1}(a_n(x_1)) \hspace{1mm} \sum_{k_2=0}^n p_{n,k_2}(a_n(x_2) \leq \varepsilon \hspace{1mm} e^{2 \mu}.
\end{align*} Moreover, we also have that
\begin{align*}
    I_{1,2} & \leq e^{2\mu}\hspace{1mm}2\Vert f \Vert_{\infty} \sum_{\left|\frac{k_1}{n}-x_1 \right|\leq \frac{\gamma}{\sqrt{2}}} p_{n,k_1}(a_n(x_1))\sum_{\left|\frac{k_2}{n}-x_2 \right|> \frac{\gamma}{\sqrt{2}}} p_{n,k_2}(a_n(x_2)) \\[0.5em]
    & \leq e^{2 \mu}\hspace{1mm}2\Vert f \Vert_{\infty}\sum_{\left|\frac{k_2}{n}-x_2 \right|> \frac{\gamma}{\sqrt{2}}}p_{n,k_2}(a_n(x_2)).
\end{align*} 
Now, we can write what follows:
\begin{align}
&   \sum_{\left|\frac{k_i}{n}-x_i \right|> \delta }p_{n,k_i}(a_n(x_i))\ =\ \sum_{\left|\frac{k_i}{n}-x_i \right|> \delta }  { \left|\frac{k_i}{n}-x_i \right|   \over   \left|\frac{k_i}{n}-x_i \right| }  p_{n,k_i}(a_n(x_i)) \notag
\\
& \leq\ \delta^{-1} \sum_{\left|\frac{k_i}{n}-x_i \right|> \delta }  \left|\frac{k_i}{n}-x_i \right|    p_{n,k_i}(a_n(x_i))  \notag
\end{align}
\begin{align}
& \leq \delta^{-1} \sum_{ \left|\frac{k_i}{n}-x_i \right|   > \delta } \left|\frac{k_i}{n}-a_n(x_i) \right|\,    p_{n,k_i}(a_n(x_i))\  \notag
\\
& +\ \delta^{-1} \left|a_n(x_i)-x_i \right|\, \sum_{ \left|\frac{k_i}{n}-x_i \right| > \delta }      p_{n,k_i}(a_n(x_i))   \notag
\\
 \label{ultima-stima}
& \leq \delta^{-1}\, {1 \over 2 \sqrt n}\ +\ \delta^{-1}\ \max_{x \in [0,1]}  \left|a_n(x)-x \right|,
\end{align}
$\delta>0$, $i=1,2$, where in the above computations we used a well-known inequality (see, e.g., \cite{AltoCam}) concerning the first order moment of the functions $p_{n,k_i}$, i.e., 
\begin{equation} \label{fundamental3}
    \sum_{k_i=0}^n \left| \frac{k_i}{n}-x \right|  p_{n,k_i}(x) \le {1\over 2\sqrt n}, \quad x \in [0,1].
\end{equation}
Observing that 
\begin{align} 
\gamma_n\ :&=\ \max_{x \in [0,1]}  \left|a_n(x)-x \right|\notag \\ 
& \label{gamma-n} =\ { n \over \mu} \ln \left[ \left( e^{\mu / n} - 1 \right) {n \over \mu} \right] - {\left( e^{\mu / n} - 1 \right) {n \over \mu } - 1 \over e^{\mu / n} - 1 }   \longrightarrow 0, \quad {\rm as} \quad n \to +\infty,  
\end{align}
and using (\ref{ultima-stima}), then we finally obtain that
$$
I_{1,1} \leq e^{2 \mu}\hspace{1mm}4\Vert f \Vert_{\infty}\hspace{1mm} \frac{\varepsilon}{\gamma^2},
$$
for $n$ sufficiently large. Since the reasoning is similar for $I_{1,3}$ and $I_{1,4}$, the thesis immediately follows. 
\end{proof}

For each $\mathbf{x} \in \mathbb{R}^d$, we consider the function $\exp_{\mu,\mathbf{x}}$ defined for every $\mathbf{t}=(t_1,\dots,t_d)$$\in \mathbb{R}^d$ by
\begin{align*} \displaystyle
    \exp_{\mu,\mathbf{x}}(\mathbf{t}) := \exp_{\mu}(\mathbf{t}) - \exp_\mu(\mathbf{x}) = e^{\mu \sum_{i=1}^d t_i}-e^{\mu \sum_{i=1}^d x_i}.
\end{align*} 
Now, using Lemma \ref{LemD} and the expressions in (\ref{exp2}) and (\ref{Conn}), we can explicitly compute the operator $\widetilde{\mathscr{G}}_n$ for the square of $\exp_{\mu,\mathbf{x}}$, for every $\mathbf{x} \in Q_d$:
\begin{align*} %\label{exp2D} 
\displaystyle
    \widetilde{\mathscr{G}}_n(\exp^2_{\mu,\mathbf{x}},\mathbf{x}) & = \widetilde{\mathscr{G}}_n(\exp^2_\mu, \mathbf{x})-2\exp_\mu(\mathbf{x})\widetilde{\mathscr{G}}_n(\exp_\mu,\mathbf{x})+\exp^2_\mu(\mathbf{x})\widetilde{\mathscr{G}}_n(\exp_0, \mathbf{x}) \notag \\[0.5em]
    & = \exp^2_\mu(\mathbf{x})-2\exp_\mu(\mathbf{x})\exp_\mu(\mathbf{x})
    +\exp^2_\mu(\mathbf{x})\widetilde{\mathscr{G}}_n(\exp_0, \mathbf{x}) \notag \\[0.5em]
    & = \exp^2_\mu(\mathbf{x})-2\exp^2_\mu(\mathbf{x})+\exp^2_\mu(\mathbf{x})\widetilde{\mathscr{G}}_n(\exp_0, \mathbf{x}) \notag \\[0.5em] & = \exp^2_\mu(\mathbf{x})\left(\widetilde{\mathscr{G}}_n(e_0, \mathbf{x})-1\right)
    \notag \\ & = e^{2\mu \sum_{i=1}^d x_i}\left(\left( \hspace{1mm} \prod_{i=1}^d \mathscr{G}_n(e_0, x_i)\right)-1\right).   
\end{align*}
In order to obtain a quantitative estimate for the approximation error of the operator $\widetilde{\mathscr{G}}_n$, we recall the definition of the modulus of continuity in the multidimensional case. For $f \in C(Q_d)$, the modulus of continuity of $f$ is given by
\begin{align*} %\label{ModMult} 
\displaystyle
    \omega(f,\delta):= \sup_{\substack{\Vert \mathbf{t}-\mathbf{x} \Vert_2 \leq \delta \\ \vphantom{\int} \mathbf{x},\mathbf{t}\in Q_d }}\left|f(\mathbf{t})-f(\mathbf{x})\right|,
\end{align*} where $\delta > 0$.

By means of a different approach, it is possible to achieve the following quantitative estimate for the order of approximation. 

\begin{theorem} \label{Teoultimo}
For every $f\in C([0,1])$, $n\in \mathbb{N}$, there holds:
\begin{align*} \displaystyle
\| \widetilde{\mathscr{G}}_n f  - f \|_{\infty} \leq & \hspace{2mm} 
\left(1+{d\over 2}\right) e^{\mu d}\, \omega\left(f,{1\over \sqrt n}\right) +   e^{\mu d}\, \Vert f\Vert_\infty{\mu^2 \over n} + e^{\mu d} \left( 1+ \mu d \right) \omega\left(f, {1\over n} \right)
\end{align*} 
\end{theorem}

\begin{proof}
We can write what follows:
%  Proceeding as in the first part of the proof of Theorem \ref{TeoUltimo} it is possible to obtain (\ref{momento}). Now, we can write what follows: 
    \begin{align*} \displaystyle
& \big|\widetilde{\mathscr{G}}_n(f,\mathbf{x})-f(\mathbf{x})\big| 
    \notag\\[0.7em] & = \hspace{1mm} \big| \widetilde{\mathscr{G}}_n(f,\mathbf{x})-\widetilde{\mathscr{G}}_n(e_0,\mathbf{x})f(\mathbf{x})+\widetilde{\mathscr{G}}_n(e_0,\mathbf{x})f(\mathbf{x})-f(\mathbf{x})\big| 
    \notag\\[0.7em] & \leq \hspace{1mm} \big|\widetilde{\mathscr{G}}_n (f,\mathbf{x})-f(\mathbf{x}) \widetilde{\mathscr{G}}_n (e_0,\mathbf{x})\big|+\big|f(\mathbf{x})\big|\big|\widetilde{\mathscr{G}}_n(e_0,\mathbf{x})-1\big| \notag \\[0.5em] 
&\leq  \sum_{k_1=0}^n {\dots} \sum_{k_d=0}^n \left|f\left(\frac{k_1}{n},\dots,\frac{k_d}{n}\right)-f(\mathbf{x})\right| e^{-\mu k_1/n} \dots e^{-\mu k_d/n} \hspace{1mm}e^{\mu x_1} \dots e^{\mu x_d}
    \notag\\[0.5em] &  \hspace{14mm} p_{n,k_1}(a_n(x_1))\dots p_{n,k_d}(a_n(x_d)) +\Vert f\Vert_\infty \big|\widetilde{\mathscr{G}}_n(e_0,\mathbf{x})-1\big| =: J_1+J_2.
\end{align*} 
Now, the estimate for $J_2$ follows immediately since, using (\ref{voro-e0}), we can obtain:
$$
J_2\ \leq\  \Vert f\Vert_\infty {\mu^2 \over n},
$$
for $n$ sufficiently large. About $J_1$, using the subadditivity of the modulus of continuity, we have 
 \begin{align*} 
 J_1 &\le  \sum_{k_1=0}^n {\dots} \sum_{k_d=0}^n \omega\left(f, \left\| \left(\frac{k_1}{n},\dots,\frac{k_d}{n}\right)-\mathbf{x} \right\|_2\right)  e^{-\mu k_1/n} \dots e^{-\mu k_d/n} e^{\mu x_1} \dots e^{\mu x_d} \\ &
    p_{n,k_1}(a_n(x_1))\dots p_{n,k_d}(a_n(x_d)) \\
& \leq \sum_{k_1=0}^n {\dots}   \sum_{k_d=0}^n \left\{ \omega\left(f, \left\| \left(\frac{k_1}{n},\dots,\frac{k_d}{n}\right)- \left(a_n(x_1),\dots,a_n(x_d)\right) \right\|_2\right) \right. \\
 & +\ \left. \omega\left(f, \left\| \left(a_n(x_1),\dots,a_n(x_d)\right)  -\mathbf{x} \right\|_2\right) \right\}  \\
& e^{-\mu k_1/n} \dots e^{-\mu k_d/n} e^{\mu x_1} \dots e^{\mu x_d} p_{n,k_1}(a_n(x_1))\dots p_{n,k_d}(a_n(x_d)) =: J_{1,A} + J_{1,B}.
\end{align*} 
Now, applying the following well-known property\hspace{1mm}
\begin{align*} \displaystyle
    \omega(f,\lambda \delta) \leq (1+\lambda)\hspace{1mm}\omega(f,\delta), \hspace{8mm} \lambda > 0,
\end{align*}
we get:
\begin{align*} 
J_{1,A}    &\le \omega\left(f, {1\over \sqrt n}\right) \sum_{k_1=0}^n {\dots} \sum_{k_d=0}^n  \left(1+\sqrt n \left\| \left(\frac{k_1}{n},\dots,\frac{k_d}{n}\right)- \left(a_n(x_1),\dots,a_n(x_d)\right) \right\|_2 \right)  \\ 
&   e^{-\mu k_1/n} \dots e^{-\mu k_d/n} e^{\mu x_1} \dots e^{\mu x_d}\  p_{n,k_1}(a_n(x_1))\dots p_{n,k_d}(a_n(x_d)) \\
    &\le e^{\mu d} \omega\left(f, {1\over \sqrt n}\right)  \left[ 1+ \sum_{k_1=0}^n {\dots} \sum_{k_d=0}^n  \sqrt n \left\| \left(\frac{k_1}{n},\dots,\frac{k_d}{n}\right)-  \left(a_n(x_1),\dots,a_n(x_d)\right)   \right\|_2 \right.\\
    & \left. p_{n,k_1}(a_n(x_1))\dots p_{n,k_d}(a_n(x_d)) \right]  \\
    & \le e^{\mu d} \omega\left(f, {1\over \sqrt n}\right)  \left[ 1+ \sum_{k_1=0}^n {\dots} \sum_{k_d=0}^n  \sqrt n  \left\{ \sum_{j=1}^d \left| \frac{k_j}{n}-a_n(x_j) \right| \right\} \right.
    \\
    & \left. p_{n,k_1}(a_n(x_1))\dots p_{n,k_d}(a_n(x_d)) \right] \\
    & \le e^{\mu d} \omega\left(f, {1\over \sqrt n}\right) \left[ 1+ \sum_{j=1}^d \left\{  \left( \prod_{i\neq j} \sum_{k_i=0}^n p_{n,k_i}(a_n(x_i)) \right) \sqrt n  \right. \right. \\
    & \left. \left. \sum_{k_j=0}^n \left| \frac{k_j}{n}-a_n(x_j) \right| p_{n,k_j}(a_n(x_j)) \right\} \right].
\end{align*} 
Now, using (\ref{fundamental3}) we can finally write
 \begin{align*} 
 J_{1,A} &\le e^{\mu d} \omega\left(f, {1\over \sqrt n}\right)  \left[ 1+ {1\over 2} \sum_{j=1}^d \prod_{i\neq j} \sum_{k_i=0}^n p_{n,k_i}(a_n(x_i)) \right] \le  \left(1+{d\over 2}\right) e^{\mu d} \omega\left(f, {1\over \sqrt n}\right).
\end{align*} 
Concerning $J_{1,B}$, we notice that (possibly with the help of some mathematical software) $\lim_{n\to +\infty} n\gamma_n={\mu \over 8}$, where $\gamma_n$ is defined in (\ref{gamma-n}). Therefore, for example, we can deduce that
$\gamma_n \le {\mu \over n}$ for $n$ large enough, and so, for every $\mathbf{x}\in Q_d$,
$\left\| \left(a_n(x_1),\dots,a_n(x_d)\right)  -\mathbf{x} \right\|_2 \le  \sum_{j=1}^d | a_n(x_j)-x_j | \le d \gamma_n \le {d \mu \over n}$. From this we deduce that
\begin{align*}
    \omega(f,\left\| \left(a_n(x_1),\dots,a_n(x_d)\right)  -\mathbf{x} \right\|_2) &\le \omega\left(f, {d\mu \over n}\right) \le \left(1+d\mu\right) \omega\left(f,{1\over n}\right),
    \end{align*}
and so, for  sufficiently large $n\in\mathbb N$,
\begin{align*}
J_{1,B} &\le \left(1+d\mu\right) \omega\left(f,{1\over n}\right)
\sum_{k_1=0}^n {\dots} \sum_{k_d=0}^n e^{-\mu k_1/n} \dots e^{-\mu k_d/n} \hspace{1mm}e^{\mu x_1} \dots e^{\mu x_d} \\ & 
p_{n,k_1}(a_n(x_1)) {\dots} p_{n,k_d}(a_n(x_d))\\
&\le e^{\mu d} \left(1+d\mu\right) \omega\left(f,{1\over n}\right)
\sum_{k_1=0}^n  
p_{n,k_1}(a_n(x_1)) {\dots}\sum_{k_d=0}^n p_{n,k_d}(a_n(x_d))\\
&= e^{\mu d} \left(1+d\mu\right) \omega\left(f,{1\over n}\right).
\end{align*}
This completes the proof.
\end{proof}

Now, in order to obtain a qualitative estimate for the approximation error of the operators $\widetilde{\mathscr{G}}_n$, we first recall the definition of the Lipschitz class. In particular, for $\alpha \in (0,1]$, the Lipschitz class of $f \in C(Q_d)$ is given by
\begin{align*} \displaystyle
   Lip(\alpha) = \{ f \in C(Q_d) \hspace{2mm} | \hspace{2mm} \omega(f, \delta) = \mathcal{O}(\delta^{\alpha}), \hspace{2mm} \delta \rightarrow 0^+\}.
\end{align*}
Furthermore, using Theorem \ref{Teoultimo}, we can derive the following desired estimate: 
\begin{align*} \displaystyle
    \Vert \widetilde{\mathscr{G}}_n f-f\Vert_\infty = \mathcal{O}(1/n^{\alpha/2} ), \hspace{4mm} n \rightarrow +\infty,
\end{align*} for any $f$ belonging to $ Lip(\alpha)$.

\section*{Acknowledgments}

{\small The first and the second authors are members of the Gruppo Nazionale per l'Analisi Matematica, la Probabilit\`a e le loro Applicazioni (GNAMPA) of the Istituto Nazionale di Alta Matematica (INdAM), of the network RITA (Research ITalian network on Approximation), and of the UMI (Unione Matematica Italiana) group T.A.A. (Teoria dell'Approssimazione e Applicazioni). 
}

\section*{Funding}

{\small The authors L. Angeloni and D. Costarelli have been partially supported within the projects: (1) 2023 GNAMPA-INdAM Project ``Approssimazione costruttiva e astratta mediante operatori di tipo sampling e loro applicazioni" CUP E53C22001930001, (2) 2024 GNAMPA-INdAM Project ``Tecniche di approssimazione in spazi funzionali con applicazioni a problemi di diffusione", CUP E53C23001670001, and (3) PRIN 2022 PNRR: ``RETINA: REmote sensing daTa INversion with multivariate functional modeling for essential climAte variables characterization", funded by the European Union under the Italian National Recovery and Resilience Plan (NRRP) of NextGenerationEU, under the Italian Ministry of University and Research (Project Code: P20229SH29, CUP: J53D23015950001). L. Angeloni has been also partially supported within the project PRIN 2022 ``EXPANSION - EXPlainable AI through high eNergy physicS for medical Imaging in ONcology'', funded by the Italian Ministry of University and Research (CUP:  J53D23002530006).

}

\section*{Conflict of interest/Competing interests}

{\small The authors declare that they have no conflict of interest and competing interest.}

\section*{Availability of data and material and Code availability}

{ \small Not applicable.}

%\nocite{*} % Includi tutte le voci nel file .bib nel documento
\bibliography{main_v3} % Includi il file .bib

\end{document}